\documentclass{article}
\usepackage{amsmath,amssymb,indentfirst}
\usepackage[latin1]{inputenc}
\usepackage{hyperref}
\usepackage{url}
\newtheorem{theorem}{Theorem}
\newtheorem{lemma}{Lemma}

\newtheorem{exe}{Example}
\newtheorem{rmk}{Remark}

\newcommand{\qed}{\hfill\mbox{\raggedright $\Box$}\medskip}

\begin{document}

\title{THE EXPONENTIAL MATRIX:\\ 
AN EXPLICIT FORMULA BY AN ELEMENTARY METHOD}

\author{Oswaldo Rio Branco de Oliveira}
\date{}
\maketitle

\begin{abstract} We show an explicit formula, with a quite easy deduction, for the exponential matrix $e^{tA}$ of a real square matrix $A$ of order $n\times n$. The elementary method developed requires neither Jordan canonical form, nor eigenvectors, nor resolution of linear systems of differential equations, nor resolution of linear systems with constant coefficients,  nor matrix inversion, nor complex integration, nor functional analysis. The basic tools  are power series and the method of partial fraction decomposition.  Two examples are given. A proof of one well-known stability result is given. 
\end{abstract}

\vspace{0,2 cm}

\hspace{- 0,6 cm}{\sl Mathematics Subject Classification: 15-01, 34-01, 34A30, 15A16, 65F60 }

\hspace{- 0,6 cm}{\sl Key words and phrases:} Ordinary Differential Equations, Linear Equations and Systems, Exponential Matrix, Functions of Matrices.

\vspace{0,3 cm}

\section{Introduction.}

\vspace{0,2 cm}

In this article we give an explicit formula, with a quite easy proof, for the exponential matrix $e^{tA}$ of a square real matrix $A$ of order $n\times n$, where $t$ is an arbitrary real number. The method developed in what folows requires neither Jordan canonical form, nor eigenvectors, nor resolution of linear systems of differential equations, nor resolution of linear systems with constant coefficients,  nor matrix inversion, nor functional analysis (as required for the integration of functions of one complex variable taking values in the Banach space of the square complex matrices of order $n\times n$). The basic tools employed in this work are (1) some basic results on complex power series and (2) the method of partial fraction decomposition.

As is well-known, the question of computing the exponential matrix $e^{tA}$ arises from the problem of finding a solution $x:\mathbb R\to \mathbb R^n$ of the constant coefficients linear system of ordinary differential  equations 
$$ \left\{\begin{array}{ll}
x'(t)=Ax(t)\\
x(0)=x_0,
\end{array}\right.$$
where $A$ is a real square matrix of order $n\times n$ and $x_0$ is a fixed point in $\mathbb R^n$. As is also well-known, the unique solution is $x(t)=e^{tA}x_0$.
 
Two of the best methods of finding the exponential matrix $e^{tA}$ are the method that employs Jordan canonical form, see Gantmacher [3, pp. 149--152], and Putzer's method, see Apostol  [1, pp. 205--208]. Some authors, short of employing Jordan canonical form, resort to the linear algebra primary decomposition theorem, see Taylor [9, pp. 146--157]. It is worth to point out that the nice method developed by Putzer requires solving another linear system of differential equations (although a handicap, this new system is not a hard one). 

Among others ways of computing $e^{tA}$ we mention Kirchner \cite{Kirchner}. In it Kirchner also find an explicit formula for $e^{tA}$. However, his method requires to compute the inverse of a matrix and this can be troublesome. On the other hand, the method provided in this article avoids matrix inversion. 

It is important to notice that a fairly sophisticated generalization of the method developed in this work is basically presented in some other texts. Such general and quite abstract method employs the Cauchy integral formula  for functions of one complex variable taking values in complex Banach spaces (e.g.,  Rudin [8, pp. 258--267]).
Nevertheless, three comments are worthwhile regarding this generalization  presented by Rudin (and some other texts). First, the {\sl Symbolic Calculus} very nicely explained in his book requires knowledge of the complex integration theory and a bit of functional analysis.  
Second, the result presented by Rudin is not explicit on how to use the method of the partial fraction decomposition in order to obtain the exponential matrix $e^{tA}$ (although that is is not hard to figure out). Third and foremost, this quite advanced approach is unnecessary in order to find the exponential matrix $e^{tA}$. 

For those who are also interested on numerical analysis and computational algorithms to evaluate the matrix $e^{tA}$, we refer Moler and Van Loan \cite{MolerVanLoan}. In it, they focus specially the cases where $A$ is a matrix of order $n\times n$ with  $n\leq 100$.


\section{Preliminaries - Power series.}

Let us denote by $z$ the complex variable in $\mathbb C$. We indicate the usual norm of a complex number $z$ by $|z|$.  Given $a_0,a_1,a_2,\ldots$ a sequence of complex numbers, it is well-known that if the real series $\sum_{n=0}^{+\infty}|a_n|$ converges then the complex series $\sum_{n=0}^{+\infty}a_n$ also does. We say that the series $\sum_{n=0}^{+\infty} a_n$ converges absolutely  if the series $\sum_{n=0}^{+\infty}|a_n|$ converges.

We also write $\sum_{n=0}^{+\infty}a_n<\infty$ if the series converges.

Let us denote a complex sequence $a_0,a_1,a_2,\ldots$ by $(a_n)=(a_0,a_1,\ldots)$, and a complex power series with complex coefficients $(a_n)$ by
$$f(z)=a_0 + a_1z + a_2z^2 + \cdots.$$

We say that such power series converges at a given point $\zeta\in \mathbb C$ if the numerical series $\sum_{n=0}^{+\infty}a_n \zeta^n$ converges in $\mathbb C$. Then we define 
$f(\zeta)=\sum_{n=0}^{+\infty} a_n\zeta^n$.


\newpage

The three results on power series shown in this section, all three enunciated for the entire complex plane, have analogous versions and proofs  that are valid for complex  power series  defined on the open ball $B(0,1)=\{z \in \mathbb C: |z|<1\}$.

\vspace{0,2 cm}

 

\begin{lemma}{\bf (Convergence and absolute convergence).} \label{L1} The complex power series
$$f(z)=\sum_{n=0}^{+\infty}a_nz^n, \textrm{where}\ z \in \mathbb C,$$
converges for all $z\in \mathbb C$ if and only if $\sum_{n=0}^{+\infty}|a_nz^n|$ converges for all $z \in \mathbb C$.
\end{lemma}
{\bf Proof.} We already know that absolute convergence implies convergence.

Let us show the other implication. The case $z=0$ is obvious. So, let us fix a point $z\neq 0$. Then we take the positive number $2|z|$. By hypothesis we have
$$\sum_{n=0}^{+\infty} a_n(2|z|)^n<\infty.$$
Thus, we have the convergence $a_n2^n|z|^n\xrightarrow{}0$ if $n\to \infty$ and then there exists $N=N(z)$ such that we have
$ |a_n||z|^n2^n\leq 1 \ \textrm{for all}\ n\geq N$. Then, it is clear that
$$\sum_{n=0}^{+\infty} |a_nz^n| \leq (|a_0|+\cdots +|a_Nz^N|)+ \sum_{n=N+1}^{+\infty}\frac{1}{2^n}<\infty.$$
$\qed$

\vspace{0,2 cm}

\begin{lemma} {\bf (Differentiation of power series).} \label{L2} Let us consider the two complex power series
$$f(z) = \sum_{n=0}^{+\infty} a_nz^n\ \ \textrm{and}\ \ g(z) =\sum_{n=1}^{+\infty}na_nz^{n-1}.$$
Then, one power series converges at every point in the complex plane if and only if the other one also does. In such cases we have 
\[f'(z) = g(z) \ \textrm{for all}\ z.\]
\end{lemma}
{\bf Proof.} Let us split the proof into two parts: convergence and differentiation.
\begin{itemize}
\item[$\diamond$] {\sf Convergence}. Evidently, both power series converge at the origin. Moreover, given any $z\in \mathbb C$ we have the inequality
$$\sum_{n=0}^{+\infty} |a_nz^n|\leq |a_0| + |z|\sum_{n=1}^{+\infty} |na_nz^{n-1}|.$$
\newpage
By the other hand, and noticing that we have $2^n\geq n$ for all  $n$, given any $z\neq 0$ we find the inequality
$$\sum_{n=1}^{+\infty} |na_nz^{n-1}|=\frac{1}{|z|}\sum_{n=1}^{+\infty}\frac{n}{2^n}|a_n(2z)^n|\leq \frac{1}{|z|}\sum_{n=0}^{+\infty} |a_n(2z)^n|.$$
Thus, by lemma \ref{L1} we conclude that $\sum_{n=0}^{+\infty} a_nz^n$ converges everywhere (in the complex plane) if and only if $\sum_{n=1}^{+\infty} na_nz^{n-1}$ converges everywhere.

\item[$\diamond$] {\sf Differentiation.} Let us suppose that both power series converge everywhere. Then, let us fix a point 
 $z\in \mathbb C$ and  $R>|z|$. Next, we consider an increment $h\in \mathbb C$ such that $0<|h|<r=R-|z|$. Fixing any $n\geq 2$ we find that
$$\left\{\begin{array}{lll}
\frac{(z+h)^n -z^n}{h}\,  =\, nz^{n-1} \, +\, h\sum\limits_{p=2}^n\binom{n}{p}z^{n-p}h^{p-2}\\
\ \textrm{and}\ \\
\left|\frac{(z+h)^n -z^n}{h}\ -\ nz^{n-1}\right|   \leq \ \frac{|h|}{r^2}\ \sum\limits_{p=2}^n\,\binom{n}{p}\,|z|^{n-p}\,r^p\  \leq \frac{|h|\,}{r^2}R^n.\\
\end{array}
\right.$$
We notice that$|z+h|<R$. Therefore we conclude that
\[\left|\,\sum_{n=0}^{+\infty} a_n\frac{(z+h)^n - z^n}{h}\,-\, \sum_{n=0}^{+\infty} na_nz^{n-1}\,\right| \ \leq \ \frac{|h|}{r^2}\sum_{n=0}^{+\infty}|a_n|R^n.\] 
Letting $h\to 0$ we find $f'(z)=g(z)$.$\qed$
\end{itemize}

As a convention, the null polynomial has degree  $-\infty$.
Let us assume  the fundamental theorem of algebra (for a proof of it, see Oliveira \cite{Oliveira2}).


\begin{lemma}{\bf (Dividing power series by polynomials).} \label{L3} Let us consider a complex and everywhere convergent  power series
$$f(z)=\sum_{n=0}^{+\infty} a_nz^n \ \textrm{where}\ z \in \mathbb C, $$ 
and a complex polynomial $p(z)=(z-\lambda_1)^{m_1}\cdots(z-\lambda_m)^{m_m}$ with distinct complex zeros $\lambda_1,\ldots,\lambda_m$ and $\textrm{degree}(p)=m_1+\cdots + m_m=n\geq 1$. Then we have
$$f(z) = q(z)p(z) + r(z)\ \textrm{for all}\ z \in \mathbb C,$$
with $q$ a power series centered at the origin that converges everywhere  in $\mathbb C$ and $r$ a complex polynomial satisfying $\textrm{degree}(r ) <\textrm{degree}(p)$. Such  $q$ and $r$ are unique.
\end{lemma}
{\bf Proof.} Let us split it into three parts: initial case, induction, and uniqueness.
\begin{itemize}
\item[$\diamond$] {\sf Initial case.} Fixing a point $z$ and $\alpha=\lambda_1$ we  factorize
$$f(z) - f(\alpha)=\sum_{n=0}^{+\infty} a_n z^n - \sum_{n= 0}^{+\infty} a_n\alpha^n = \sum_{n=1}^{+\infty}a_n(z^n - \alpha^n)\ \ \ \ \ \ \ \ \  \ \ \ \ \ \ \ \ \ \ $$
$$\ \ \ \ \ \ \ \ \ \ \ \ \ \ = (z-\alpha)\sum_{n=1}^{+\infty}a_n(z^{n-1} + z^{n-2}\alpha + \cdots + z\alpha^{n-2} + \alpha^{n-1}).$$

\newpage

Now, let us fix $\rho$ satisfying $0\leq \max\{|z|,|\alpha|\}< \rho$. Then we have
$$ \sum_{\substack{n \in \{1,\ldots, N\} \\  k\in \{0,\ldots,n-1\}}}|a_n||z|^{n-1 -k}\,|\alpha|^k\leq \sum_{n=1}^{+\infty}|a_n|n\rho^{n-1},\ \textrm{for all}\ N.$$
This last series is finite by Lemma \ref{L2}. Under such condition, in order to compute the sum $\sum_{n,k}a_nz^{n-1-k}\alpha^k$, we can freely associate its terms (see Lang \cite{Lang}, Oliveira \cite{Oliveira1}). 
Therefore, we may write
$$\sum_{n\geq 1}a_n(z^{n-1} + z^{n-2}\alpha + \cdots + z\alpha^{n-2} + \alpha^{n-1})=\sum_{n\geq 0} b_nz^n,$$
with $(b_n)$ a complex sequence. 
This is true for all $z$. This shows
$$f(z)= (z-\alpha)Q_1(z) + f(\alpha),\ \textrm{with}\ Q_1(z)= \sum_{n\geq 0} b_nz^n\ \textrm{for all}\ z.$$
\item[$\diamond$] {\sf Induction.} Iterating the initial case,  it is fairly trivial to see that we have $f(z)= q(z)p(z) + r(z)$ with $q(z)$ a power series centered at the origin and convergent over $\mathbb C$, and $r(z)$ a polynomial with degre$(r)<$degree$(p)$. 

\item[$\diamond$]{\sf Uniqueness.} Let us suppose that we have two decompositions
$$f(z)= q_1(z)p(z) + r_1(z) \ \ \ \textrm{and}\ \ \ f(z)= q_2(z)p(z) + r_2(z),$$
both satisfying the requirements in the statement. Then we find
$$(q_2 -q_1)p = r_1-r_2.$$
Let us suppose that $r_1-r_2$ is not the null polynomial. Let $\lambda$ be a zero of multiplicity $j\geq 1$ of $p$. It is not difficult to see that there exists 
$$\lim_{z \to \lambda}\frac{r_1(z) - r_2(z)}{ (z - \lambda)^j}.$$ 
 This shows that $\lambda$ is a zero of $r_1 -r_2$ with multiplicity $k\geq j$. This is true for every zero of $p$. Thus we have $\textrm{degree}(r)\geq \textrm{degree}(p)$, a contradiction. 

So, we have $r_1=r_2$ and then $q_1(z)p(z)=q_2(z)p(z)$ for all $z$. From this, it is quite clear that we have $q_1(z)=q_2(z)$ for all $z$.
$\qed$

\end{itemize}

\begin{rmk} Lemma \ref{L3} has an analogous result, valid if we suppose that $f(z)$ is a complex power series convergent in $B(0,1)$ and 
$p(z)=(z-\lambda_1)^{m_1}\cdots (z-\lambda_m)^{m_m}$ is a polynomial with its zeros $\lambda_1,\ldots,\lambda_m$ inside $B(0,1)$. It is not difficult to figure out  a proof of such result that almost replicates the one right above. 
\end{rmk}

\newpage

\section{Partial fraction decomposition.}

\begin{lemma} {\bf (Useful derivatives).}\label{L4} Let us consider  a point $\lambda\in \mathbb C$, the complex variable $z \in \mathbb C$, an open ball $B(\lambda,\rho)=\{z\in \mathbb C: |z-\lambda|<\rho\}$ centered at $\lambda$ with radius $\rho>0$, an integer $N\geq 1$,  and an infinitely differentiable function 
$$g:B(\lambda,\rho)\to \mathbb C.$$
 Then we have
$$\frac{d^k\left\{(z-\lambda)^N\right\}}{dz^k}\Big|_{z=\lambda}=\left\{\begin{array}{ll}
0\ \textrm{for all} \ k\neq N\\
\\
N!\ \textrm{if}\ k=N
\end{array}\right.$$
and
$$\frac{d^k\left\{(z-\lambda)^Ng(z)\right\}}{dz^k}\Big|_{z=\lambda}=0\ \textrm{for all}\ k=0,\ldots,N-1.$$
\end{lemma}
{\bf Proof.} It is trivial and we leave it to the reader.$\qed$

\vspace{0,2 cm}

A complex polynomial $p=p(z)$ is monic if its dominant coefficient is $1$. 

\vspace{0,2 cm}

\begin{theorem} {\bf (Partial fraction decomposition).} \label{T1} Let $f$ and $q$ be power series convergent over $\mathbb C$, and  $p$ and $r$ be polynomials,  all as in Lemma \ref{L3}. That is,
$$f(z) = q(z)p(z) + r(z)\ \textrm{with degree$(r)<n=$\ degree$(p)$}.$$
Let us suppose that $p$ is monic and degree$(p)=n\geq 1$. Let $\lambda_1,\ldots,\lambda_m$ be the distinct zeros of $p(z)$, with respective multiplicities $m_1,\ldots,m_m$. Let us write
$$p(z)= (z-\lambda_1)^{m_1}\cdots(z-\lambda_m)^{m_m}.$$
Then, there are $n$ constants $C_{11},\ldots,C_{1,m_1}, \ldots,C_{m_11},\ldots,C_{m_1m_m}$ such that we have the decomposition
$$\frac{f(z)}{p(z)} = q(z)+ 
\left[\frac{C_{11}}{z - \lambda_1}+\cdots + \frac{C_{1m_1}}{(z-\lambda_1)^{m_1}}\right]+\cdots + 
\left[\frac{C_{m1}}{z - \lambda_m}+\cdots + \frac{C_{mm_m}}{(z-\lambda_m)^{m_m}}\right], $$
for all $z\in \mathbb C\setminus \{\lambda_1,\ldots,\lambda_m\}$.
These constants are unique and given by the derivatives
$$\left\{\begin{array}{ll}
C_{jk_j}=\frac{g_j^{(m_j-k_j)}(\lambda_j)}{(m_j - k_j)!}\
\textrm{for all}\ j=1,\ldots,m\ \textrm{and all}\ k_j=1,\ldots,m_j,
\\ \textrm{where}\\
g_{j}(z)= \frac{f(z)(z-\lambda_j)^{m_j}}{p(z)}= \frac{f(z)}{\prod\limits_{l\neq j}(z-\lambda_l)^{m_l}}.
\end{array}
\right. $$
\end{theorem}


\hspace{-0,6 cm}{\bf Proof.} Let us split it into three parts: decomposition,  formula, and uniqueness.
\begin{itemize}
\item[$\diamond$] {\sf Decomposition.} Let us show it by induction on $n$, where $n=\textrm{degree}(p)$. The case $n=1$ is trivial.

Given $n\geq 2$, let us suppose that the decomposition holds for polynomials with degree smaller or equal to $n-1$. Then, given $p=p(z)$ with degree equal to $n$, by the fundamental theorem of algebra we may write 
$$\left\{\begin{array}{ll}
p(z)=(z-\lambda_1)^{m_1}P(z), \\
\textrm{with}\ 0\leq \textrm{degree}(P) <n \ \textrm{and}\ P(\lambda_1)\neq 0.
\end{array}\right.$$
 Next, we may write
$$r(z)= \frac{r(\lambda_1)}{P(\lambda_1)}P(z)+ \left[r(z)- \frac{r(\lambda_1)}{P(\lambda_1)}P(z)\right].$$
Clearly, the polynomial inside brackets has a zero at $z=\lambda_1$.

Thus, the Euclidean algorithm guarantees  a polynomial $s(z)$ satisfying
$$r(z) -\frac{r(\lambda_1)}{P(\lambda_1)}P(z)= (z - \lambda_1)s(z)\  
\textrm{with}\ 
\textrm{degree$(s)\leq n-2$}.$$

Hence, by the identity  $p(z)=(z-\lambda_1)^{m_1}P(z)$ we arrive at
$$\frac{r(z)}{p(z)}=\frac{r(\lambda_1)/P(\lambda_1)}{(z-\lambda_1)^{m_1}}+\frac{s(z)}{(z - \lambda_1)^{m_1-1}P(z)}.$$
We have $\textrm{degree}(s)\leq n-2<n-1=\textrm{degree}[ (z - \lambda_1)^{m_1-1}P(z)]$. Therefore, arguing by induction we arrive at the 
desired decomposition.

\item[$\diamond$]{\sf Formula.} Let us fix $j=1$. Given $k_1\in\{1,\ldots, m_1\}$, following the above decomposition and the definition of $g_1$ (given in the statement) we find
$$g_1(z)=C_{11}(z-\lambda_1)^{m_1 -1} +\cdots +C_{1k_1}(z-\lambda_1)^{m_1 -k_1}+\cdots + 
C_{1m_1}$$
$$\ \ \ \ + (z-\lambda_1)^{m_1}\left[q(z)+\sum_{\substack{2\,\leq\, j\,\leq m\\ 1\,\leq\, k_j\,\leq\, m_j}}\frac{C_j{k_j}}{(z-\lambda_j)^{m_j}}\right].$$
We notice that the function inside brackets is infinitely differentiable on a neighborhood of the point $z=\lambda_1$.

From Lemma \ref{L4}  it follows that
$$g_1^{(m_1 -k_1)}(\lambda _1)=(m_1 - k_1)!C_{1k_1}.$$
The argument for $j=2,\ldots,m$ is analogous.

\item[$\diamond$] {\sf Uniqueness.} The uniqueness of  $C_{11},\ldots, C_{mm_m}$ follows from the formula.$\qed$
\end{itemize}

\newpage

\section{ The explicit formula for $e^{tA}$.}


The following result is well-known and we omit its proof 
(see Apostol \cite{Apostol}). 

\begin{lemma}{\bf(Cayley-Hamilton Theorem).}  \label {T2} Let $A$ be a real square matrix of order $n\times n$ and $p_A(z)=\det(zI - A)=z^n + a_{n-1}z^{n-1} +  \cdots + a_1z + a_0$ its characteristic polynomial. Then we have
$$p_A(A)= A^n +a_{n-1}A^{n-1}+ \cdots + a_1A+a_0I=0.$$
\end{lemma}

\begin{theorem}{\bf (Explicit Formula for $e^{tA}$).} \label{T3} Let $A$ be a real square matrix of order $n\times n$ with characteristic polynomial 
$$p_A(z)=(z-\lambda_1)^{m_1}\cdots (z -\lambda_m)^{m_m},$$
where $\lambda_1,\ldots,\lambda_m$ are the distinct zeros of $p_A$ with $m_1,\ldots,m_m$ their respective algebraic multiplicities. For each $j=1,\ldots,m$ and  each $k_j=1,\ldots, m_j$, let us consider the polynomial (a total of $n$ polynomials)
$$\ \ \ \ \ \ \ \ \ \ \ \ \ \ \ p_{jk_j}(z)= (z-\lambda_j)^{m_j - k_j}\prod_{l\neq j}(z - \lambda_l)^{m_l} \ \ \  \left[=\frac{p_A(z)}{(z-\lambda_j)^{k_j}}\right].$$
Then we have (to simplify, we omit the set where the indices take values) 
$$e^{tA}= \sum C_{jk_j}p_{jk_j}(A),$$
where 
$$C_{jk_j}=\frac{1}{(m_j - k_j)!}\frac{d^{m_j - k_j}}{dz^{m_j - k_j}}\left\{ \frac{e^{tz}(z-\lambda_j)^{m_j}}{p_A(z)}\right\}\Big|_{z=\lambda_j}. $$
\end{theorem}
{\bf Proof.} Fixing an arbitrary $t\in \mathbb R$, the function $z\mapsto e^{tz}$ is given by a power series that converges over the entire complex plane. By Lemma \ref{L3} we have
$$\ \ e^{tz}=q(z)p_A(z) + r(z),\ \ \textrm{with}\ \left\{\begin{array}{ll}
q\ \textrm{a power series that converges over}\ \mathbb C,\\
r\ \textrm{a polynomial with}\ \textrm{degree}(r)<\ \textrm{degree}(p_A).
\end{array}\right.
$$
Since $A$ comutes with powers of $A$ and with the identity matrix, we find that
$$e^{tA}=q(A)p_A(A)+r(A).$$
The Cayley-Hamilton theorem shows that $p_A(A)=0$. Thus we have
$$e^{tA}=r(A).$$
By the decomposition in Theorem \ref{T1} we may write
$$\frac{r(z)}{p_A(z)}=\sum_{\substack{1\leq j\leq m\\ 1 \leq k_j\leq m_j}} \frac{C_{j k_j}}{(z-\lambda_j)^{k_j}}.$$
Thus we are allowed to conclude that
$$r(z)= \sum C_{j k_j}p_{jk_j}(z)\ \ \textrm{and}\ \ r(A)=\sum C_{jk_j}p_{jk_j}(A).$$
The proof is complete.$\qed$

\newpage

\section{Examples.}

\begin{exe} -  \label{E1}  Let us compute $e^{tA}$ for
$$A=\left(\begin{array}{ll}
1 & 2\\
4 & 3\end{array}
\right).$$
\end{exe}
{\bf Solution.} The characteristic polynomial is
$$p(z) =\left|\begin{array}{ll}
z - 1 & \  - 2\\
\  - 4 & z -3\end{array}
\right|=(z-1)(z-3) -8=(z+1)(z-5). $$
Following Theorem \ref{T1}, and its notation, we write
$$\frac{e^{tz}}{(z+1)(z-5)}= q(z)+\frac{r(z)}{(z+1)(z-5)}= q(z) + \frac{\alpha}{z+1} + \frac{\beta}{z-5},$$
where $q=q(z)$ is complex power series centered at the origin and $\alpha$ and $\beta$ are complex constants. These constants are given by
$$\alpha = -\frac{e^{-t}}{6} \ \ \ \textrm{and}\ \ \beta=\frac{e^{5t}}{6}.$$
Thus, we arrive at
$$e^{tA}= \frac{e^{5t}}{6}(A +I) -\frac{e^{-t}}{6}(A-5I)=  \frac{e^{5t}}{6}\left(\begin{array}{ll}
2 & 2\\
4 & 4\end{array}
\right) -\frac{e^{-t}}{6}\left(\begin{array}{ll}
-4 & \ \, 2\\
\ \, 4 & -2\end{array}
\right).$$ 

\begin{exe} - {\bf Real matrices of order $3\times 3$.} \label{E2}  Let us compute $e^{tA}$ when $A$ is a real matrix of order $3\times 3$, with characteristic polynomial $p(z)$. 
\end{exe}
{\bf Solution.} Let us start with the case $p(z)=(z-\lambda)^2(z-\mu)$ where $\lambda\neq \mu$. Following Theorem \ref{T1}, and its notation, we may write
$$\frac{e^{tz}}{(z-\lambda)^2(z-\mu)}=q(z) + \frac{\alpha}{z-\mu} + \frac{\beta}{z-\lambda} + \frac{\gamma}{(z-\lambda)^2},$$
where $q=q(z)$ is a complex power series centered at the origin and $\alpha$, $\beta$, and $\gamma$ are complex constants. These constants are given by
$$\alpha = \frac{e^{\mu t}}{(\mu - \lambda)^2},\ \ \  \gamma= \frac{e^{\lambda t}}{\lambda - \mu},$$
and
$$\beta = \frac{d}{dz}\left\{\frac{e^{tz}}{z-\mu}\right\}\Big|_{z=\lambda}= \frac{te^{tz}(z-\mu) - e^{tz}}{(z - \mu)^2}\Big|_{z=\lambda}=\frac{e^{\lambda t}[t(\lambda -\mu) -1]}{(\lambda - \mu)^2}.$$ 
The exponential matrix searched for is
$$e^{tA} = \frac{e^{\mu t}}{(\mu - \lambda)^2}(A-\lambda I)^2 + \frac{e^{\lambda t}[t(\lambda -\mu) -1]}{(\lambda - \mu)^2}(A-\lambda I)(A-\mu I) + \frac{e^{\lambda t}}{\lambda - \mu}(A-\mu I).$$
%

We leave the case when $p(z)=(z -\lambda)(z - \mu)(z - \nu)$ has three distinct  roots  to the reader. We also leave the case when $p(z)=(z-\lambda)^3$ has a root of multiplicity three to the reader.


\begin{exe} -  {\bf A Stability Result.} \label{E3}  Let $A$ be a $n\times n$ real matrix and $p(z)$ be its characteristic polynomial. Let us suppose that all the characteristic roots $\lambda_1,\ldots,\lambda_m$ of $p(z)$ have negative real  part. Then, let us show that there exist a constant $\alpha <0$ and a constant $C>0$ such that the unique solution of the linear system
$$\left\{\begin{array}{ll}
x'(t)=Ax(t)\\ 
x(0)=x_0,\\
\end{array}\right. $$
satisfy
$$|x(t)|\leq Ce^{\alpha t}|x_0|,\ \textrm{for all}\ t \in \mathbb R \ \textrm{and all}\ x_0\in \mathbb R^n.$$
\end{exe}
{\bf Proof.} From the hypotheses it follows that there exists $\alpha <0$ satisfying
$$ \textrm{Re}(\lambda_j) < \alpha, \ \textrm{for all}\ j =1,\ldots,m.$$
By employing Theorem \label{T3} and its notation we have
$$e^{tA}= \sum C_{jk_j}p_{jk_j}(A).$$
We notice that the matrix $p_{jk_j}(A)$ does not depend on the variable $t$. Hence, there exists $M>0$ such that we have
$$|p_{jk_j}(A)|\leq M,\ \textrm{for all allowed}\ j  \ \textrm{and}\ k_j.$$
Thus we have
$$|e^{tA}| \leq M\sum |C_{jk_j}|.$$
The coefficient $C_{jk_j}$ depends on the variable $t$ and it is given by 
$$ C_{jk_j}=\frac{1}{(m_j - k_j)!}\frac{d^{(m_j - k_j)}}{dz^{m_j - k_j}}\left\{ \frac{e^{tz}(z-\lambda_j)^{m_j}}{p_A(z)}\right\}\Big|_{z=\lambda_j}. $$
Thus, there exists a complex constant $D_{jk_j}$ such that we have
$$C_{jk_j}=D_{jk_j}e^{\lambda_j t}.$$ 
Evidently, there exists $N>0$ such that we have
$$|D_{jk_j}|\leq N,\ \textrm{for all allowed}\ j  \ \textrm{and}\ k_j.$$
Hence, we arrive at
$$|e^{tA}|\leq MN\sum e^{\textrm{Re}(\lambda_j)t}\leq MN\sum e^{\alpha t} \leq Ce^{\alpha t},\ \textrm{for some}\ C>0.$$
Therefore, we conclude that the solution $x(t)=e^{tA}x_0$ satisfy
$$|x(t)|\leq Ce^{\alpha t}|x_0|,\ \textrm{for all} \ t \in \mathbb R\ \textrm{and all}\ x_0\in \mathbb R^n.$$
$ \qed$

\newpage

\paragraph{Acknowledgments.}  The author is thankful to Professor G. Terra for his comments and for reference \cite{Rudin}.

\bigskip

\noindent\textit{Departamento de Matemática,
Universidade de São Paulo\\
Rua do Matão 1010 - CEP 05508-090\\
São Paulo, SP - Brasil\\
oliveira@ime.usp.br}

\bigskip

\end{document}